\documentclass[11pt]{article}
\usepackage{latexsym,amssymb,amsmath,amsfonts,amsthm}
\usepackage{graphics,graphicx,mathrsfs,subfigure}
\usepackage{color,xspace}
\newcommand{\R}{{\mat R}}

\newcommand{\Sp}{{\mat S}}

\newcommand{\be}{\begin{eqnarray}}
\newcommand{\ben}{\begin{eqnarray*}}
\newcommand{\en}{\end{eqnarray}}
\newcommand{\enn}{\end{eqnarray*}}

\newcommand{\mat}{\mathbb}

\newtheorem{theorem}{Theorem}[section]
\newtheorem{lemma}[theorem]{Lemma}

\newtheorem{proposition}[theorem]{Proposition}

\definecolor{rot}{rgb}{1,0,0}
\definecolor{hw}{rgb}{0,0,1}

\begin{document}
\renewcommand{\theequation}{\arabic{section}.\arabic{equation}}
%\begin{titlepage}
\title{\bf
Determination of the potential by a fixed angle scattering data
}
\author{ Suliang Si\thanks{School of Mathematics and Statistics, Shandong University of Technology,
Zibo, 255000, China ({\tt sisuliang@amss.ac.cn})}}
\date{}
%\end{titlepage}

%\vspace{.2in}

%\begin{document}
%\renewcommand{\theequation}{\arabic{section}.\arabic{equation}}

\maketitle

\begin{abstract}
In this paper, we study the fixed angle scattering problems raised by Rakesh
and Salo (SIAM J Math Anal 52(6):5467–5499, 2020) and (Inverse Probl
36(3):035005, 2020). We show that a compactly
supported potential is uniquely determined by the far field pattern  at a fixed angle, which is an open problem. Our method is based on a novel Carleman estimate
and  the ideas introduced by Bukhgeim and Klibanov on the use of
Carleman estimates for inverse problems.
\end{abstract}

%\begin{keywords}
%Uniqueness, inverse scattering problem, far-field pattern, penetrable obstacle,
%transmission problem, interior transmission problem, embedded obstacles
%\end{keywords}
%
%\begin{AMS}
%78A46, 65C30
%\end{AMS}

%\pagestyle{myheadings}
%\thispagestyle{plain}
%\markboth{J. Yang, B. Zhang, and H. Zhang}{Determining penetrable obstacles with embedded objects}

\section{Introduction}

We are mainly concerned with the inverse scattering problem of determining a compactly supported
unknown potential in the Schr{\"o}dinger equation in \(\R^n\), \(n\geq 2\).
Assume
the incident field is given by the  plane wave
\[\hat{u}^i(x,k,\omega) = e^{ik \omega \cdot x},\]
where $k > 0$ is a frequency and $\omega \in {\Sp}^{n-1}$ the direction of propagation. Then the scattering problem is to the full field \(\hat{u}\) such that 
\begin{equation}\label{0u}
-\Delta \hat{u}-k^2\hat{u}+V(x)\hat{u}=0, \, \,x\in\R^n,
\end{equation}
\begin{equation}\label{1u}
\hat{u}(x,k,d)=e^{ik\omega\cdot x}+\hat{u}^{s}(x,k,\omega),
\end{equation}
\begin{equation}\label{r}
r^{\frac{n-1}{2}}(\partial_r \hat{u}^{s}-ik\hat{u}^{s})=0, \ \,r:=|x|\rightarrow +\infty,
\end{equation}
where \(V\in C_c^\infty(\R^n)\) and (\ref{r}) is \textit{the Sommerfeld radiation condition} which
guarantees that the scattered wave \(\hat{u}^s\) is outgoing.
 The well-posedness of the forward scattering problem (\ref{0u})-(\ref{r}) can be conveniently found in \cite{CK19}.
Writing $x = r\theta$ where $r \geq 0$ and $\theta \in {\Sp}^{n-1}$, the scattered wave \(\hat{u}^s\) has the asymptotics

\[
\hat{u}^s(r\theta,k,\omega) = e^{ik r} r^{-\frac{n-1}{2}} \hat{u}_{\infty}(\theta, k, \omega; V) + o(r^{-\frac{n-1}{2}}) \quad \text{ as } \,  r \to \infty.
\]
The function $\hat{u}_{\infty}(\theta, k, \omega; V)$ is called the \textit{far field pattern}, corresponding to the potential $V$. One could interpret $\hat{u}_{\infty}(\theta, k, \omega; V)$ as a scattering measurement for $V$ that corresponds to sending a plane wave at frequency $k > 0$ propagating in the direction $\omega \in {\Sp}^{n-1}$ and measuring the scattered wave in the direction $\theta \in {\Sp}^{n-1}$.
We show that a compactly supported potential \(V\) is uniquely determined by the far field pattern at a fixed angle \(\omega\).

\begin{theorem}\label{T0}
Fix $\omega \in {\Sp}^{n-1}$, $n \geq 2$, and let $V_1, V_2 \in C_c^\infty(\mathbb{R}^n)$ be real valued. If 
\[
\hat{u}_{\infty}(\theta, k, \omega; V_1) = \hat{u}_{\infty}(\theta, k, \omega; V_2) 
\]
for all $k > 0$ and $\theta \in {\Sp}^{n-1}$, then $V_1 = V_2$.
\end{theorem}

In inverse scattering problems the objective is to determine certain properties of a scatterer from
measurements that are made far away.
We formulate four fundamental inverse scattering problems, related to recovering a potential from (partial) knowledge of its quantum mechanical scattering amplitude:

\begin{enumerate}
\item \textbf{Full data.} Recover \( V \) from \( \hat{u}_{\infty}(\theta, k, \omega; V) \)  for all \(\theta \in {\Sp}^{n-1}\), \( k > 0 \), \( \omega \in {\Sp}^{n-1} \).
\item \textbf{Fixed frequency.} Recover \( V \) from \( \hat{u}_{\infty}(\theta, k, \omega; V) \) with \( k > 0 \) fixed.
\item \textbf{Backscattering.} Recover \( V \) from \( \hat{u}_{\infty}(\omega, k, \omega; V) \) for all \( k > 0 \), \( \omega \in {\Sp}^{n-1} \).
\item \textbf{Fixed angle.} Recover \( V \) from \( \hat{u}_{\infty}(\theta, k, \omega; V) \) where \( \omega \in {\Sp}^{n-1} \) is fixed.
\end{enumerate}
The full data problem is formally overdetermined when \( n \geq 2 \), since one seeks to recover a function of \( n \) variables from a function of \( 2n - 1 \) variables. Similarly, the fixed frequency problem is formally overdetermined when \( n \geq 3 \) (it is formally determined when \( n = 2 \)). Both of these problems have been solved; we only mention that one can determine \( V \) from the high frequency asymptotics of \( \hat{u}_{\infty} \) \cite{Sa82} and that the fixed frequency problem is equivalent to a variant of the inverse conductivity problem of Calderón addressed in \cite{Bu08, SU87}. There have been many related works and we refer to \cite{Uh92, No08, Uh14} for references.

The backscattering and the fixed angle inverse scattering problems are formally determined in any dimension (both the unknown and the data depend on \( n \) variables). The one-dimensional case is well understood \cite{Ma11, DT79}. Known results for \( n \geq 2 \) include uniqueness for potentials that are small or belong to a generic set \cite{ER92, St92, MU08, B+20}, recovery of main singularities \cite{GU93, OPS01, Ru01}, identification of the zero potential in fixed angle scattering \cite{BLM89}, and recovery of angularly controlled potentials from backscattering data \cite{RU14}. See the references in \cite{RU14, Me18} for further results. However, these problems remain open in general.

For inverse coefficient problems of the hyperbolic equation, concerning the uniqueness and
 the stability, Bukhgeim and Klibanov \cite{BK81} proposed a fundamental method based on
 what is called a Carleman estimate. Thus Carleman estimates became a fundamental
 tool for establishing uniqueness and stability for inverse problems \cite{BY17,IY01}.
Rakesh and M. Salo showed that a compactly supported potential is
uniquely determined by the far field pattern generated by plane waves coming
from exactly two opposite directions \cite{RM20}. These results are proved using Carleman estimates
and adapting the ideas introduced by Bukhgeim and Klibanov on the use of
Carleman estimates for inverse problems. Later, they extended the methods of \cite{RM20} to establish an equivalence between the
frequency domain and the time domain formulations of the problem. So they obtained that that a compactly
supported potential is uniquely determined by its scattering amplitude for two opposite fixed angles.  In this work, we establish a novel Carleman estimate. Based on this inequality, we  prove that a compactly
supported potential is uniquely determined by its scattering amplitude for only one fixed angle, which solves an open problem.

The paper is organized as follows. In section \ref{time}, we state the time
domain setting for the fixed angle scattering problem. Section \ref{SEC} is an intoduction to a novel Carleman estimate. In Section, we prove Theorem \ref{T2}. Finally, Appendix  contains the derivation of a new Carleman estimate (\ref{CAR}) without inter information.

\section{The time domain setting}\label{time}
Assume that the support of $V(x)$ is contained in \(B:=\{x\in\R^n|\,\, |x|<1\}\).
We consider the initial value problem with a plane wave source:
\begin{equation}\label{U}
\partial_t^2 U-\Delta U + VU = 0, \quad (x, t) \in \mathbb{R}^n \times \mathbb{R}, 
\end{equation}
\begin{equation}\label{U0}
U(x, t) = \delta(t - x \cdot \omega), \quad x \in \mathbb{R}^n, \, \, t \ll 0,
\end{equation}
where \(\delta\) is the Dirac function.
This was studied in \cite{RU14} and the following proposition \ref{pro} is a consequence of the arguments in the proof of  \cite[Theorem 1]{RU14}.

\begin{proposition} \label{pro}
The initial value problem (\ref{U}), (\ref{U0}) has a unique distributional solution \(U(x, t, \omega)\) given by
\begin{equation}
U(x, t, \omega) = \delta(t - x \cdot \omega) + u(x, t, \omega)H(t - x \cdot \omega),
\end{equation}
where \(u(x, t, \omega)\), a smooth function on the region \(t \geq x \cdot \omega\), is the unique solution of the characteristic initial value problem:
\begin{equation}\label{u}
\partial_t^2 u -\Delta u+ Vu = 0, \quad (x, t) \in \mathbb{R}^n \times \mathbb{R}, \, t > x \cdot \omega, 
\end{equation}
\begin{equation}\label{u1}
u(x, x \cdot \omega, \omega) = -\frac{1}{2} \int_{-\infty}^0 V(x + \sigma \omega) d\sigma, \quad x \in \mathbb{R}^n, 
\end{equation}
\begin{equation}\label{u2}
u(x, t, \omega) = 0, \quad x \in \mathbb{R}^n, \, x \cdot \omega < t \ll 0. 
\end{equation}
Also, for any real \(T\), on the region \(\{(x, t) : x \cdot \omega \leq t \leq T\}\), \(|u(x, t, \omega)|\) is bounded above by a continuous function of \(\|V\|_{C^{n+4}}\).
\end{proposition}
From the above proposition, we know there exists a constant \(M>0\) satisfying
\begin{equation}
|u(x, t,\omega)|\leq C\|V\|_{C^{n+4}}\leq M \quad \text{for all } x\cdot \omega\leq t\leq T.
\end{equation}
Let \(Q_+=\{(x, t) |\, x\in B, \, \,  x \cdot \omega \leq t \leq T\}\) and 
\(\Sigma_+=\{(x, t) |\, x\in \partial B, \, \,  x \cdot \omega \leq t \leq T\}\).
\begin{theorem}[One-plane wave data]\label{T2}
Let $u_i(x,t,\omega)$ be the solutions of (\ref{u})-(\ref{u2}) with $V=V_i$, $i=1,2$. If
\[
u_1(x, t, \omega)  =  u_2(x, t, \omega), \quad (x, t) \in (\partial B \times \mathbb{R}) \cap \{ t \geq x \cdot \omega \},
\]
then
\[V_1=V_2 \quad  \text{for all }x\in\R^n. \]
\end{theorem}
This is our main  result, the proof is placed in section \ref{Proof}. 
The following theorem shows that the scattering amplitude for a fixed direction \(\omega\in {\Sp}^{n-1}\) and
the boundary measurements in the wave equation problem  are equivalent information, which can be found in \cite{RM20}.
\begin{theorem}
Let \( n \geq 2 \) and fix \( \omega \in {\Sp}^{n-1} \), \(k_0>0\). For any real valued \( V_1, V_2 \in C_c^\infty(\mathbb{R}^n) \) with support in \( \bar{B} \), let $u_i(x,t,\omega)$ be the solutions of (\ref{u})-(\ref{u2}) with $V=V_i$, $i=1,2$, one has
\[
\hat{u}_{\infty}(\theta, k, \omega; V_1) = \hat{u}_{\infty}(\theta, k, \omega; V_2) \quad  \text{ for all }  \,  k \geq k_0 \text{ and }  \, \theta \in {\Sp}^{n-1}
\]
if and only if
\[
u_1(x, t, \omega)  =  u_2(x, t, \omega) \quad \text{ for } \quad (x, t) \in (\partial B \times \mathbb{R}) \cap \{ t \geq x \cdot \omega \}.
\]
\end{theorem}
From the above two theorems, we deduce  Theorem \ref{T0}  immediately.

\section{A new Carleman estimate for hyperbolic equation}\label{SEC}

Throughout this article, \(M> 0\) denote generic constants which are independent of parameter \(s\).
For convenience, we use \(\{(x,t)|\,\, x\in B, \, t=T\}\) by \(\{t=T\}\) and \(\{(x,t)| \,\, x\in B,\, t=x\cdot \omega\}\) by \(\{t=x\cdot \omega\}\) respectively.
 In deriving Carleman estimates, we will use the following simple integration by parts results on \( Q_+ \) and other sets having a similar form. If \( v \) is smooth in \( Q_+ \), then
\[
\int_{Q_+} \partial_t v \, \mathrm{d}x \, \mathrm{d}t = \int_{\{t = T\}} v \, \mathrm{d}S - \frac{1}{\sqrt{2}} \int_{\{t=x\cdot \omega\}} v \, \mathrm{d}S,
\]
and if \( \textbf{v} \) is a smooth vector field on \( Q_+ \) with values in \( \mathbb{R}^n \), then (with \( \nabla \) denoting the gradient in \( x \) variables)
\[
\int_{Q_+} \nabla \cdot \textbf{v} \, \mathrm{d}x \, \mathrm{d}t = \int_{\Sigma_+} \textbf{v} \cdot \nu \, \mathrm{d}S + \frac{1}{\sqrt{2}} \int_{\{t=x\cdot \omega\}} \textbf{v} \cdot \omega \, \mathrm{d}S.
\]
Let \(Q_+^\infty=\{(x, t) |\, x\in B, \, \,  t\geq x \cdot \omega \}\) and 
\(\Sigma_+^\infty=\{(x, t) |\, x\in \partial B, \, \, t\geq x \cdot \omega \}\).

Assume \(x_0\notin \overline{B}\) and \(0<t_0<T\). Let \(\beta\in(0,1)\) and define 
\begin{equation}
\psi (x,t)=|x-x_0|^2-\beta (t-t_0)^2+C_0\quad \text{and for}\quad \lambda>0, \quad \varphi (x,t)=e^{\lambda\psi (x,t)}, 
\end{equation}
where \(C_0>0\) is chosen such that \(\psi\geq 1\) in \(Q_+^\infty\).
\begin{lemma}\label{Car}
Assume \(t_0\) is sufficiently large. Then there exists \(s_0\) such that for all \(s>s_0\), we have
\begin{equation}\label{Car1}
\begin{split}
&s^{1/2}\int_{t=x\cdot\omega} e^{2s\varphi}|\partial_tv+\nabla v\cdot\omega|^2  dS +s^{5/2}\lambda^2\int_{t=x\cdot\omega} e^{2s\varphi }\varphi^2(|\partial_t\psi|^2-|\nabla\psi|^2)|v|^2dS\\
& + s  \int_{Q_+^\infty} e^{2s\varphi} \left( |\partial_t v|^2 + |\nabla v|^2 \right) dxdt + s^3  \int_{Q_+^\infty} e^{2s\varphi} |v|^2 dxdt \\
&
\leq M \int_{Q_+^\infty}e^{2s\varphi}|\partial_t^2v-\Delta v+Vv|^2dxdt
\end{split}
\end{equation}
for all \(v\in H^2(Q_+^\infty)\) satisfying \(v=\partial_tv=0\) on \(\Sigma_+^\infty\).
\end{lemma}
It is sufficient to prove Lemma \ref{Car} in the case where \(V\equiv 0\). Indeed, we assume that we already established
the inequality (\ref{Car1}).
Since \(V\in C_c^\infty(\mathbb{R}^n)\), we have
\begin{equation}
|\partial_t^2v-\Delta v|^2\leq |\partial_t^2v-\Delta v+Vv-Vv|^2\leq 2 |\partial_t^2v-\Delta v+Vv|^2+M|v|^2.
\end{equation}
By choosing \(s\) large, we can absorb the term
\[ \int_{Q_+} e^{2s\varphi} |v|^2 dxdt\]
into the left-hand side of (\ref{Car1}) with \(V=0\).
Let \(0<t_0<T\).
In order to prove the Carleman estimates (\ref{Car1}), we set \[z=e^{s\varphi}v \quad \mbox{for all } (x,t)\in Q_+.\] 
Then, we introduce the conjugate operator $P$ defined by
\begin{equation}
Pz = e^{s\varphi} (\partial_t^2-\Delta) \left( e^{-s\varphi} z \right).
\end{equation}
Some easy computations give
\begin{equation}
\begin{aligned}
Pz &= \partial_t^2 z - 2s\lambda\varphi\left( \partial_t z \partial_t\psi - \nabla z \cdot \nabla\psi \right) + s^2\lambda^2\varphi^2 z\left( |\partial_t\psi|^2 - |\nabla\psi|^2 \right) - \Delta z \\
&\quad - s\lambda\varphi z\left( \partial_t^2\psi - \Delta\psi \right) - s\lambda^2\varphi z\left( |\partial_t\psi|^2 - |\nabla\psi|^2 \right)-s\partial_t z \\
&= P_1 z + P_2 z + R_1z,
\end{aligned}
\end{equation}
where
\begin{align}
P_1 z=\partial_t^2 z- \Delta z+ s^2\lambda^2\varphi^2 z\left( |\partial_t\psi|^2 - |\nabla\psi|^2 \right),
\end{align}
\begin{align}
P_2 z &= (\alpha - 1)s\lambda\varphi z(\partial_t^2\psi - \Delta\psi) - s\lambda^2\varphi z(|\partial_t\psi|^2 - |\nabla\psi|^2) \notag \\
&\quad - 2s\lambda\varphi(\partial_t z\partial_t\psi - \nabla z\cdot\nabla\psi)
\end{align}
and
\begin{align}
R_1 z &= -\alpha s\lambda\varphi z(\partial_t^2\psi - \Delta\psi).
\end{align}
Let
\begin{equation}
\alpha \in \left( \frac{2\beta}{\beta + n}, \frac{2}{\beta + n} \right).
\end{equation}
Since we have
\begin{equation}
\int_{Q_+} \left( |P_1 z|^2 + |P_2 z|^2 \right) dxdt + 2 \int_{Q_+} P_1 z P_2 z dxdt = \int_{Q_+} |P z - R_1 z|^2 dxdt, 
\end{equation}
the main part of the proof is then to bound from below the cross-term
\[
\int_{Q_+} P_1 z P_2 z dxdt=\sum_{i,k=1}^3 I_{i,k}.
\]
We calculate the six terms \(I_{i,k}\), \(i,k=1,2,3\) by integrating by parts with respect to  \((x,t)\).

Integrations by part in time give easily
\[
\begin{aligned}
I_{11} &= \int_{Q_+} \partial_t^2 z \left( (\alpha - 1)s\lambda\varphi z (\partial_t^2\psi - \Delta\psi) \right) dxdt \\
&=(\alpha - 1)s\lambda\int_{t=T}\partial_t z z\varphi (\partial_t^2\psi - \Delta\psi)dS-\frac{(\alpha - 1)s}{\sqrt{2}}\lambda\int_{t=x\cdot\omega}\partial_t z z\varphi (\partial_t^2\psi - \Delta\psi)dS
\\
&-(\alpha - 1)s\lambda\int_{Q_+} |\partial_tz|^2\varphi (\partial_t^2\psi - \Delta\psi)dxdt\\
&-\frac{(\alpha - 1)s\lambda}{2}\int_{t=T} |z|^2\partial_t\big(\varphi (\partial_t^2\psi - \Delta\psi)\big)dS+\frac{(\alpha - 1)s\lambda}{2\sqrt{2}}\int_{t=x\cdot\omega} |z|^2\partial_t\big(\varphi (\partial_t^2\psi - \Delta\psi)\big)dS
\\
&+\frac{(\alpha - 1)s\lambda^2}{2}\int_{Q_+} |z|^2\varphi\partial^2_t\psi (\partial_t^2\psi - \Delta\psi)dxdt+\frac{(\alpha - 1)s\lambda^3}{2}\int_{Q_+} |z|^2|\varphi\partial_t\psi|^2 (\partial_t^2\psi - \Delta\psi)dxdt
\end{aligned}
\]
Similarly, one has
\[
\begin{aligned}
I_{12} &= \int_{Q_+} \partial_t^2 z \left( -s\lambda^2\varphi z \left( |\partial_t\psi|^2 - |\nabla\psi|^2 \right) \right) dxdt \\
&=-s\lambda^2\int_{t=T} \partial_t z  z\varphi(|\partial_t\psi|^2 - |\nabla\psi|^2)dS+\frac{s\lambda^2}{\sqrt{2}}\int_{t=x\cdot\omega} \partial_t z  z\varphi(|\partial_t\psi|^2 - |\nabla\psi|^2)dS
\\
&+s\lambda^2\int_{Q_+} |\partial_tz|^2 \left(\varphi(|\partial_t\psi|^2 - |\nabla\psi|^2)\right)dxdt\\
& +\frac{s\lambda^2}{2}\int_{t=T} |z|^2 \partial_t\left(\varphi(|\partial_t\psi|^2 - |\nabla\psi|^2)\right)dS-\frac{s\lambda^2}{2\sqrt{2}}\int_{t=x\cdot\omega} |z|^2 \partial_t\left(\varphi(|\partial_t\psi|^2 - |\nabla\psi|^2)\right)dS
\\
&-\frac{s\lambda^2}{2}\int_{Q_+} |z|^2 \partial^2_t\left(\varphi(|\partial_t\psi|^2 - |\nabla\psi|^2)\right)dxdt\\
& - s\lambda^2 \int_{Q_+} \varphi |z|^2 |\partial_t^2\psi|^2 dxdt  - \left( 2 + \frac{1}{2} \right) s\lambda^3 \int_{Q_+} \varphi |z|^2 |\partial_t\psi|^2 \partial_t^2\psi dxdt
\\
& + \frac{s\lambda^3}{2}\int_{Q_+} \varphi |z|^2 |\nabla\psi|^2 \partial_t^2\psi dxdt  - \frac{s\lambda^4}{2}  \int_{Q_+} \varphi |z|^2 |\partial_t\psi|^2 \left( |\partial_t\psi|^2 - |\nabla\psi|^2 \right) dxdt
\end{aligned}
\]
and
\[
\begin{aligned}
I_{13} &=  \int_{Q_+} \partial_t^2 z \left( -2s\lambda\varphi \left( \partial_t z \partial_t\psi - \nabla z \cdot \nabla\psi \right) \right) dxdt \\
&= -s\lambda \int_{t=T} |\partial_tz|^2\varphi\partial_t\psi dS+\frac{s\lambda}{\sqrt{2}} \int_{t=x\cdot\omega} |\partial_tz|^2\varphi\partial_t\psi dS
\\
&+2s\lambda \int_{t=T} \partial_tz \varphi\nabla z\cdot\nabla\psi dS-\frac{2s\lambda}{\sqrt{2}} \int_{t=x\cdot\omega} \partial_tz \varphi\nabla z\cdot\nabla\psi dS
\\
&+ s\lambda \int_{Q_+} \varphi |\partial_t z|^2 \partial_t^2\psi \, dxdt + s\lambda^2  \int_{Q_+} \varphi |\partial_t z|^2 |\partial_t\psi|^2 \, dxdt \\
&\quad + s\lambda  \int_{Q_+} \varphi |\partial_t z|^2 \Delta\psi \, dxdt + s\lambda^2 \int_{Q_+} \varphi |\partial_t z|^2 |\nabla\psi|^2 \, dxdt \\
&\quad - 2s\lambda^2  \int_{Q_+} \varphi \partial_t z \partial_t\psi \nabla z \cdot \nabla\psi \, dxdt\\
&-s\lambda\int_{\Sigma_+} \varphi|\partial_t z|^2\partial_\nu \psi dS-\frac{s\lambda}{\sqrt{2}}\int_{t=x\cdot\omega} \varphi|\partial_t z|^2\nabla \psi\cdot\omega dS.
\end{aligned}
\]
Furthermore, by Green’s formula and integration by parts, we obtain
\[
\begin{aligned}
I_{21} &= \int_{Q_+} -\Delta z \left( (\alpha - 1)s\lambda\varphi z (\partial_t^2\psi - \Delta\psi) \right) dxdt \\
&= -(1 - \alpha)s\lambda \int_{Q_+} \varphi |\nabla z|^2 (\partial_t^2\psi - \Delta\psi) dxdt+ \frac{(1 - \alpha)}{2} s\lambda^2  \int_{Q_+} \varphi |z|^2 \Delta\psi (\partial_t^2\psi - \Delta\psi) dxdt \\
&\quad + \frac{(1 - \alpha)}{2} s\lambda^3  \int_{Q_+} \varphi |z|^2 |\nabla\psi|^2 (\partial_t^2\psi - \Delta\psi) dxdt\\
&+(1 - \alpha)s\lambda  \int_{\Sigma_+} \partial_\nu z\varphi z (\partial_t^2\psi - \Delta\psi) dS+\frac{(1 - \alpha)s\lambda}{\sqrt{2}}  \int_{t=x\cdot\omega} \nabla z\cdot\omega\varphi z (\partial_t^2\psi - \Delta\psi) dS\\
&
-\frac{(1 - \alpha)s\lambda}{2}\int_{\Sigma_+}|z|
^2\partial_\nu(\varphi(\partial_t^2\psi - \Delta\psi))dS-\frac{(1 - \alpha)s\lambda}{2\sqrt{2}}\int_{t=x\cdot\omega}|z|
^2\nabla(\varphi(\partial_t^2\psi - \Delta\psi))\cdot\omega dS
.
\end{aligned}
\]
On the other hand,
\[
\begin{aligned}
I_{22} &=\int_{Q_+} -\Delta z \left( -s\lambda^2\varphi z \left( |\partial_t\psi|^2 - |\nabla\psi|^2 \right) \right) dxdt \\
&= -s\lambda^2  \int_{Q_+} \varphi |\nabla z|^2 \left( |\partial_t\psi|^2 - |\nabla\psi|^2 \right) dxdt \\
&\quad - \frac{s\lambda^2}{2}  \int_{Q_+} \varphi |z|^2 \Delta\left( |\nabla\psi|^2 \right) dxdt \\
&\quad + \frac{s\lambda^3}{2}  \int_{Q_+} \varphi |z|^2 \Delta\psi \left( |\partial_t\psi|^2 - |\nabla\psi|^2 \right) dxdt \\
&\quad + \frac{s\lambda^4}{2} \int_{Q_+} \varphi |z|^2 |\nabla\psi|^2 \left( |\partial_t\psi|^2 - |\nabla\psi|^2 \right) dxdt \\
&\quad - s\lambda^3  \int_{Q_+} \varphi |z|^2 \nabla\psi \cdot \nabla\left( |\nabla\psi|^2 \right) dxdt\\
&+s\lambda^2  \int_{\Sigma_+} \partial_\nu z\varphi  z \left( |\partial_t\psi|^2 - |\nabla\psi|^2 \right) dS+\frac{s\lambda^2}{\sqrt{2}}  \int_{t=x\cdot\omega} \nabla z\cdot\omega\varphi  z \left( |\partial_t\psi|^2 - |\nabla\psi|^2 \right) dS
\\
&
-\frac{s\lambda^2}{2}  \int_{\Sigma_+} | z|^2 \partial_\nu \Big(\varphi\left( |\partial_t\psi|^2 - |\nabla\psi|^2 \right)\Big) dS-\frac{s\lambda^2}{2\sqrt{2}}  \int_{t=x\cdot\omega} | z|^2 \nabla \Big(\varphi\left( |\partial_t\psi|^2 - |\nabla\psi|^2 \right)\Big)\cdot\omega dS
.
\end{aligned}
\]
Using the fact that $z|_{\partial\Omega\times(0,T)} = 0$, $\nabla z = (\partial_\nu z)\nu$ and $|\nabla z|^2 = |\partial_\nu z|^2$ on $\Sigma_+$, we obtain
\[
\begin{aligned}
I_{23} &= \int_{Q_+} -\Delta z \left( -2s\lambda\varphi \left( \partial_t z \partial_t\psi - \nabla z \cdot \nabla\psi \right) \right) dxdt \\
&=-s\lambda \int_{t=T}|\nabla z|^2\varphi\partial_t\psi dS+\frac{s\lambda}{\sqrt{2}} \int_{t=x\cdot\omega}|\nabla z|^2\varphi\partial_t\psi dS\\
&+ s\lambda \int_{Q_+} \varphi |\nabla z|^2 (\partial_t^2\psi - \Delta\psi) dxdt + 2s\lambda^2  \int_{Q_+} \varphi |\nabla\psi \cdot \nabla z|^2 dxdt \\
&\quad - 2s\lambda^2  \int_{Q_+} \varphi \partial_t z \partial_t\psi \nabla z \cdot \nabla\psi dxdt + s\lambda^2 \int_{Q_+} \varphi |\nabla z|^2 \left( |\partial_t\psi|^2 - |\nabla\psi|^2 \right) dxdt \\
&\quad  + 4s\lambda  \int_{Q_+} \varphi |\nabla z|^2 dxdt\\
&+2s\lambda  \int_{\Sigma_+} \partial_\nu z\varphi \partial_tz\partial_t\psi dS+\frac{2s\lambda}{\sqrt{2}}  \int_{t=x\cdot\omega} \nabla z\cdot\omega\varphi \partial_tz\partial_t\psi dS
\\
&-
2s\lambda \int_{\Sigma_+} \partial_\nu z\varphi \nabla z\cdot\nabla\psi dS-
\frac{2s\lambda}{\sqrt{2}} \int_{t=x\cdot\omega} \nabla z\cdot\omega\varphi \nabla z\cdot\nabla\psi dS
\\
&+
 s\lambda \int_{\Sigma_+} \varphi |\nabla z|^2 \nabla\psi \cdot \nu \, dS+
 \frac{s\lambda}{\sqrt{2}} \int_{t=x\cdot\omega} \varphi |\nabla z|^2 \nabla\psi \cdot \omega \, dS
\end{aligned}
\]
One easily writes
\[
\begin{aligned}
I_{31} &=  \int_{Q_+} s^2\lambda^2\varphi^2 z \left( |\partial_t\psi|^2 - |\nabla\psi|^2 \right) \left( (\alpha - 1)s\lambda\varphi z (\partial_t^2\psi - \Delta\psi) \right) dxdt \\
&= (\alpha - 1)s^3\lambda^3  \int_{Q_+} \varphi^3 |z|^2 (\partial_t^2\psi - \Delta\psi) \left( |\partial_t\psi|^2 - |\nabla\psi|^2 \right) dxdt
\end{aligned}
\]
and
\[
\begin{aligned}
I_{32} &= \int_{Q_+} s^2\lambda^2\varphi^2 z \left( |\partial_t\psi|^2 - |\nabla\psi|^2 \right) \left( -s\lambda^2\varphi z \left( |\partial_t\psi|^2 - |\nabla\psi|^2 \right) \right) dxdt \\
&= -s^3\lambda^4  \int_{Q_+} \varphi^3 |z|^2 \left( |\partial_t\psi|^2 - |\nabla\psi|^2 \right)^2 dxdt.
\end{aligned}
\]
Finally, some integrations by part enable to obtain
\[
\begin{aligned}
I_{33} &=  \int_{Q_+} s^2\lambda^2\varphi^2 z \left( |\partial_t\psi|^2 - |\nabla\psi|^2 \right) \left( -2s\lambda\varphi \left( \partial_t z \partial_t\psi - \nabla z \cdot \nabla\psi \right) \right) dxdt \\
&=-s^3\lambda^3\int_{t=T} |z|^2\varphi^3 (|\partial_t\psi|^2 - |\nabla\psi|^2)\partial_t\psi dS+\frac{s^3\lambda^3}{\sqrt{2}}\int_{t=x\cdot\omega} |z|^2\varphi^3 (|\partial_t\psi|^2 - |\nabla\psi|^2)\partial_t\psi dS\\
&+ s^3\lambda^3  \int_{Q_+} \varphi^3 |z|^2 (\partial_t^2\psi - \Delta\psi) \left( |\partial_t\psi|^2 - |\nabla\psi|^2 \right) dxdt \\
&\quad + 2s^3\lambda^3  \int_{Q_+} \varphi^3 |z|^2 \left( \partial_t^2\psi |\partial_t\psi|^2 + 2|\nabla\psi|^2 \right) dxdt \\
&\quad + 3s^3\lambda^4  \int_{Q_+} \varphi^3 |z|^2 \left( |\partial_t\psi|^2 - |\nabla\psi|^2 \right)^2 dxdt\\
&+ s^3\lambda^3\int_{\Sigma_+} \varphi^3 (|\partial_t\psi|^2 - |\nabla\psi|^2)|z|^2\partial_\nu\psi dS+ \frac{s^3\lambda^3}{\sqrt{2}}\int_{t=x\cdot\omega} \varphi^3 (|\partial_t\psi|^2 - |\nabla\psi|^2)|z|^2\nabla\psi\cdot\omega dS.
\end{aligned}
\]
Gathering all the terms that have been computed, we get
\begin{equation}
\begin{aligned}
\int_{Q_+} P_1 z P_2 z \, dxdt
&= 2s\lambda\int_{Q_+} \varphi |\partial_t z|^2 \partial_t^2\psi \, dxdt - \alpha s\lambda \int_{Q_+} \varphi |\partial_t z|^2 (\partial_t^2\psi - \Delta\psi) dxdt \\
&\quad + 2s\lambda^2 \int_{Q_+} \varphi \left( |\partial_t z|^2 |\partial_t\psi|^2 - 2\partial_t z \partial_t\psi \nabla z \cdot \nabla\psi + |\nabla\psi \cdot \nabla z|^2 \right) dxdt \\
&\quad + 4s\lambda \int_{Q_+} \varphi |\nabla z|^2 dxdt + \alpha s\lambda \int_{Q_+} \varphi |\nabla z|^2 (\partial_t^2\psi - \Delta\psi) dxdt \\
&\quad + 2s^3\lambda^4 \int_{Q_+} \varphi^3 |z|^2 \left( |\partial_t\psi|^2 - |\nabla\psi|^2 \right)^2 dxdt \\
&\quad + 2s^3\lambda^3 \int_{Q_+} \varphi^3 |z|^2 \left( \partial_t^2\psi |\partial_t\psi|^2 + 2|\nabla\psi|^2 \right) dxdt \\
&\quad + \alpha s^3\lambda^3 \int_{Q_+} \varphi^3 |z|^2 (\partial_t^2\psi - \Delta\psi) \left( |\partial_t\psi|^2 - |\nabla\psi|^2 \right) dxdt+X_1
\\
& +X_{t=T}+X_{t=x\cdot\omega}+X_{\Sigma_+}\\
\end{aligned}
\end{equation}
where $X_1$ gathers the non-dominating terms and satisfies
\[
|X_1| \leq M s \lambda^4 \int_{Q_+} \varphi |z|^2 \, dx dt.
\]
The quantities \(X_{t=T}\), \(X_{t=x\cdot\omega}\) and \(X_{\Sigma_+}\) correspond to boundary terms from integration by parts at  \(\{t=T\}\), \(\{t=x\cdot \omega\}\), and \(\Sigma_+\), with their explicit forms presented below.

Since \(\psi=|x-x_0|^2-\beta(t-t_0)^2+C_0\) and \(
\alpha \in \left( \frac{2\beta}{\beta + n}, \frac{2}{\beta + n} \right)\), we have 
\begin{equation}
\begin{aligned}
2s\lambda \int_{Q_+} \varphi |\partial_t z|^2 \partial_t^2 \psi \,dxdt
- \alpha s\lambda \int_{Q_+} \varphi |\partial_t z|^2 \big(\partial_t^2 \psi - \Delta \psi\big) \,dxdt \\
\quad + 4s\lambda \int_{Q_+} \varphi |\nabla z|^2 \,dxdt
+ \alpha s\lambda \int_{Q_+} \varphi |\nabla z|^2 \big(\partial_t^2 \psi - \Delta \psi\big) \,dxdt \\
\ge M s\lambda \int_{Q_+} \varphi |\partial_t z|^2 \,dxdt
+ M s\lambda \int_{Q_+} \varphi |\nabla z|^2 \,dxdt.
\end{aligned}
\end{equation}
On the one hand, about the first order derivative terms, one can notice that
\begin{equation}
\begin{aligned}
2s\lambda^2 \int_{Q_+} \varphi \big(|\partial_t z|^2 |\partial_t \psi|^2 - 2\partial_t z \partial_t \psi \nabla z\cdot\nabla \psi + |\nabla \psi\cdot\nabla z|^2\big) dxdt \\
= 2s\lambda^2 \int_{Q_+} \varphi \big(\partial_t z \partial_t \psi - \nabla z\cdot\nabla \psi\big)^2 dxdt \ge 0.
\end{aligned}
\end{equation}
Considering the $0$-th order terms, we observe that
\begin{equation}
\begin{aligned}
&2s^3\lambda^4 \int_{-T}^{T}\int_{\Omega} \varphi^3 |w|^2 \big(|\partial_t\psi|^2 - |\nabla\psi|^2\big)^2 dxdt + \alpha s^3\lambda^3 \int_{-T}^{T}\int_{\Omega} \varphi^3 |w|^2 \big(\partial_t^2\psi - \Delta\psi\big)\big(|\partial_t\psi|^2 - |\nabla\psi|^2\big) dxdt \\
&\quad + 2s^3\lambda^3 \int_{-T}^{T}\int_{\Omega} \varphi^3 |w|^2 \big(\partial_t^2\psi|\partial_t\psi|^2 + 2|\nabla\psi|^2\big) dxdt+X_1 \\
&\geq Ms^3\lambda^3 \int_{-T}^{T}\int_{\Omega} \varphi^3 |w|^2 F_{\lambda,t_0}(\phi) dxdt,
\end{aligned}
\end{equation}
where
\begin{equation}
\begin{aligned}
F_{\lambda,t_0}(\phi) &= 2\lambda\big(|\partial_t\psi|^2 - |\nabla\psi|^2\big)^2 + 2\big(\partial_t^2\psi|\partial_t\psi|^2 + 2|\nabla\psi|^2\big) + \alpha\big(\partial_t^2\psi - \Delta\psi\big)\big(|\partial_t\psi|^2 - |\nabla\psi|^2\big) \\
&= 2\lambda\big(|\partial_t\psi|^2 - |\nabla\psi|^2\big)^2 + \big(2\partial_t^2\psi + \alpha(\partial_t^2\psi - \Delta\psi)\big)\big(|\partial_t\psi|^2 - |\nabla\psi|^2\big) + 2\big(\partial_t^2\psi + 2\big)|\nabla\psi|^2 \\
&= 2\lambda X^2 + \big(2\partial_t^2\psi + \alpha(\partial_t^2\psi - \Delta\psi)\big)X + 16(1-\beta)|x-x_0|^2.
\end{aligned}
\end{equation}
with $X = |\partial_t\psi|^2 - |\nabla\psi|^2$.

We now give the explicit formula for \(X_{t=x\cdot\omega}\) as
\begin{equation}\label{txw}
\begin{split}
X_{t=x\cdot\omega}&=\frac{(1-\alpha )s}{\sqrt{2}}\lambda\int_{t=x\cdot\omega}(\partial_t z+\nabla z\cdot\omega) z\varphi (\partial_t^2\psi - \Delta\psi)dS
\\
&+\frac{(\alpha - 1)s\lambda}{2\sqrt{2}}\int_{t=x\cdot\omega} |z|^2\Big(\partial_t\big(\varphi (\partial_t^2\psi - \Delta\psi)\big)+\nabla\big(\varphi (\partial_t^2\psi - \Delta\psi)\big)\cdot\omega\Big)dS
\\
&+\frac{s\lambda^2}{\sqrt{2}}\int_{t=x\cdot\omega} (\partial_t z +\nabla z\cdot\omega) z\varphi(|\partial_t\psi|^2 - |\nabla\psi|^2)dS
\\
&-\frac{s\lambda^2}{2\sqrt{2}}\int_{t=x\cdot\omega} |z|^2 \Big(\partial_t\left(\varphi(|\partial_t\psi|^2 - |\nabla\psi|^2)\right)+\nabla\left(\varphi(|\partial_t\psi|^2 - |\nabla\psi|^2)\right)\cdot\omega\Big) dS
\\
&
+\frac{s\lambda}{\sqrt{2}} \int_{t=x\cdot\omega} |\partial_tz|^2\varphi(\partial_t\psi-\nabla \psi\cdot\omega) dS+\frac{2s\lambda}{\sqrt{2}} \int_{t=x\cdot\omega} \partial_tz\partial_t\psi\varphi\nabla z\cdot\omega dS
\\
&-\frac{2s\lambda}{\sqrt{2}} \int_{t=x\cdot\omega} (\partial_tz+\nabla z\cdot\omega) \varphi\nabla z\cdot\nabla\psi dS
\\
&+\frac{s\lambda}{\sqrt{2}} \int_{t=x\cdot\omega}|\nabla z|^2\varphi(\partial_t\psi+\nabla \psi\cdot\omega) dS+\frac{s^3\lambda^3}{\sqrt{2}}\int_{t=x\cdot\omega} |z|^2\varphi^3 (|\partial_t\psi|^2 - |\nabla\psi|^2)(\partial_t\psi+\nabla \psi\cdot\omega) dS. 
\end{split}
\end{equation}
We choose \(t_0\) sufficiently large to ensure that \(\partial_t\psi+\nabla \psi\cdot\omega\) is positive at \(\{t=x\cdot\omega\}\).
A simple calculation gives
\begin{equation}
\begin{split}
&\frac{s\lambda}{\sqrt{2}} \int_{t=x\cdot\omega} |\partial_tz|^2\varphi(\partial_t\psi-\nabla \psi\cdot\omega) dS+\frac{2s\lambda}{\sqrt{2}} \int_{t=x\cdot\omega} \partial_tz\partial_t\psi\varphi\nabla z\cdot\omega dS
-\frac{2s\lambda}{\sqrt{2}} \int_{t=x\cdot\omega} (\partial_tz+\nabla z\cdot\omega) \varphi\nabla z\cdot\nabla\psi dS\\
&\quad
+\frac{s\lambda}{\sqrt{2}} \int_{t=x\cdot\omega}|\nabla z|^2\varphi(\partial_t\psi+\nabla \psi\cdot\omega) dS\\
&\geq \frac{s\lambda}{\sqrt{2}} \int_{t=x\cdot\omega} \varphi\partial_t\psi|\partial_tz+\nabla z\cdot\omega|^2  dS
-\frac{s\lambda}{\sqrt{2}} \int_{t=x\cdot\omega} |\partial_tz|^2\varphi\nabla \psi\cdot\omega dS+\frac{s\lambda}{\sqrt{2}} \int_{t=x\cdot\omega}|\nabla z\cdot\omega|^2\varphi\nabla \psi\cdot\omega dS\\
&-\frac{2s\lambda}{\sqrt{2}} \int_{t=x\cdot\omega} (\partial_tz+\nabla z\cdot\omega) \varphi\nabla z\cdot\nabla\psi dS.
\end{split}
\end{equation}
If \(\partial_tz+\nabla z\cdot\omega\equiv 0\) for all \(t=x\cdot\omega\), then 
\begin{equation}
\begin{split}
&\frac{s\lambda}{\sqrt{2}} \int_{t=x\cdot\omega} \varphi\partial_t\psi|\partial_tz+\nabla z\cdot\omega|^2  dS
-\frac{s\lambda}{\sqrt{2}} \int_{t=x\cdot\omega} |\partial_tz|^2\varphi\nabla \psi\cdot\omega dS+\frac{s\lambda}{\sqrt{2}} \int_{t=x\cdot\omega}|\nabla z\cdot\omega|^2\varphi\nabla \psi\cdot\omega dS\\
&-\frac{2s\lambda}{\sqrt{2}} \int_{t=x\cdot\omega} (\partial_tz+\nabla z\cdot\omega) \varphi\nabla z\cdot\nabla\psi dS\equiv 0.
\end{split}
\end{equation}
If \(\partial_tz+\nabla z\cdot\omega\neq 0\) for all \(t=x\cdot\omega\).
Since
\begin{equation}
\begin{split}
-\frac{2s\lambda}{\sqrt{2}} \int_{t=x\cdot\omega} (\partial_tz+\nabla z\cdot\omega) \varphi\nabla z\cdot\nabla\psi dS\geq 
-\frac{s\lambda}{\sqrt{2}} \int_{t=x\cdot\omega} \varphi|\partial_tz+\nabla z\cdot\omega |^2  dS-\frac{s\lambda}{\sqrt{2}} \int_{t=x\cdot\omega}  \varphi |\nabla\psi||\nabla z|^2 dS,
\end{split}
\end{equation}
Choosing \(t_0\) large enough, we obtain
\begin{equation}
\begin{split}
&\frac{s\lambda}{\sqrt{2}} \int_{t=x\cdot\omega} \varphi\partial_t\psi|\partial_tz+\nabla z\cdot\omega|^2  dS
-\frac{s\lambda}{\sqrt{2}} \int_{t=x\cdot\omega} |\partial_tz|^2\varphi\nabla \psi\cdot\omega dS+\frac{s\lambda}{\sqrt{2}} \int_{t=x\cdot\omega}|\nabla z|^2\varphi\nabla \psi\cdot\omega dS\\
&-\frac{2s\lambda}{\sqrt{2}} \int_{t=x\cdot\omega} (\partial_tz+\nabla z\cdot\omega) \varphi\nabla z\cdot\nabla\psi dS\geq 0.
\end{split}
\end{equation}
Therefore we always obtain
\begin{equation}
\begin{split}
&\frac{s\lambda}{\sqrt{2}} \int_{t=x\cdot\omega} |\partial_tz|^2\varphi(\partial_t\psi-\nabla \psi\cdot\omega) dS+\frac{2s\lambda}{\sqrt{2}} \int_{t=x\cdot\omega} \partial_tz\partial_t\psi\varphi\nabla z\cdot\omega dS
\\
&-\frac{2s\lambda}{\sqrt{2}} \int_{t=x\cdot\omega} (\partial_tz+\nabla z\cdot\omega) \varphi\nabla z\cdot\nabla\psi dS
+\frac{s\lambda}{\sqrt{2}} \int_{t=x\cdot\omega}|\nabla z|^2\varphi(\partial_t\psi+\nabla \psi\cdot\omega) dS\geq 0
\end{split}
\end{equation}
The first term on the right-hand side of (\ref{txw})
\begin{equation}
\frac{(1-\alpha )s}{\sqrt{2}}\lambda\int_{t=x\cdot\omega}(\partial_t z+\nabla z\cdot\omega) z\varphi (\partial_t^2\psi - \Delta\psi)dS
\end{equation}
can be controlled by
\begin{equation}
\frac{s\lambda}{\sqrt{2}}\int_{t=x\cdot\omega} \varphi\partial_t\psi|\partial_tz+\nabla z\cdot\omega|^2  dS
\end{equation}
 and
\begin{equation}
\frac{s^3\lambda^3}{\sqrt{2}}\int_{t=x\cdot\omega} |z|^2\varphi^3 (|\partial_t\psi|^2 - |\nabla\psi|^2)(\partial_t\psi+\nabla \psi\cdot\omega) dS.
\end{equation}
when \(t_0\) is sufficiently large.
The second term on the right-hand side of (\ref{txw})
\begin{equation}
\frac{(\alpha - 1)s\lambda}{2\sqrt{2}}\int_{t=x\cdot\omega} |z|^2\Big(\partial_t\big(\varphi (\partial_t^2\psi - \Delta\psi)\big)+\nabla\big(\varphi (\partial_t^2\psi - \Delta\psi)\big)\cdot\omega\Big)dS
\end{equation}
can be controlled by  
\begin{equation}
\frac{s^3\lambda^3}{\sqrt{2}}\int_{t=x\cdot\omega} |z|^2\varphi^3 (|\partial_t\psi|^2 - |\nabla\psi|^2)(\partial_t\psi+\nabla \psi\cdot\omega) dS.
\end{equation}
The third term on the right-hand side of (\ref{txw})
\begin{equation}
\frac{s\lambda^2}{\sqrt{2}}\int_{t=x\cdot\omega} (\partial_t z +\nabla z\cdot\omega) z\varphi(|\partial_t\psi|^2 - |\nabla\psi|^2)dS
\end{equation}
can be controlled by
\begin{equation}
\frac{s\lambda}{\sqrt{2}} \int_{t=x\cdot\omega} \varphi\partial_t\psi|\partial_tz+\nabla z\cdot\omega|^2  dS
\end{equation}
and 
\begin{equation}
\frac{s^3\lambda^3}{\sqrt{2}}\int_{t=x\cdot\omega} |z|^2\varphi^3 (|\partial_t\psi|^2 - |\nabla\psi|^2)(\partial_t\psi+\nabla \psi\cdot\omega) dS.
\end{equation}
In conclusion, for sufficiently large \(t_0\), we have
\[X_{t=x\cdot\omega}\geq 0.\]
The expression for \(X_{\Sigma_+}\) is
\begin{equation}
\begin{split}
X_{\Sigma_+}&=-s\lambda\int_{\Sigma_+} \varphi|\partial_t z|^2\partial_\nu \psi dS+(1 - \alpha)s\lambda  \int_{\Sigma_+} \partial_\nu z\varphi z (\partial_t^2\psi - \Delta\psi) dS
\\
&
-\frac{(1 - \alpha)s\lambda}{2}\int_{\Sigma_+}|z|
^2\partial_\nu(\varphi(\partial_t^2\psi - \Delta\psi))dS
+s\lambda^2  \int_{\Sigma_+} \partial_\nu z\varphi  z \left( |\partial_t\psi|^2 - |\nabla\psi|^2 \right) dS
\\
&
-\frac{s\lambda^2}{2}  \int_{\Sigma_+} | z|^2 \partial_\nu \Big(\varphi\left( |\partial_t\psi|^2 - |\nabla\psi|^2 \right)\Big) dS\\
&
+2s\lambda  \int_{\Sigma_+} \partial_\nu z\varphi \partial_tz\partial_t\psi dS
-
2s\lambda \int_{\Sigma_+} \partial_\nu z\varphi \nabla z\cdot\nabla\psi dS
\\
&+
 s\lambda \int_{\Sigma_+} \varphi |\nabla z|^2 \nabla\psi \cdot \nu \, dS+ s^3\lambda^3\int_{\Sigma_+} \varphi^3 (|\partial_t\psi|^2 - |\nabla\psi|^2)|z|^2\partial_\nu\psi dS.
\end{split}
\end{equation}
As \(v = \partial_tv = 0\) on \(\Sigma_+\), we get \(z = \partial_tz = 0\) on \(\Sigma_+\), hence
\[X_{\Sigma_+}=0.\]
Based on the above results and Formula (1), we have
\begin{equation}\label{P12}
\begin{split}
\int_{Q_+}P_1z P_2zdxdt\geq  M s\lambda \int_{Q_+} \varphi |\partial_t z|^2 \,dxdt
+ M s\lambda \int_{Q_+} \varphi |\nabla z|^2 \,dxdt+M s^3\lambda^3 \int_{Q_+} F_{\lambda,t_0}(\phi)\varphi^3 | z|^2 \,dxdt+X_T.
\end{split}
\end{equation}
Integrating $P_1 z \partial_t z$ by parts over $Q_+$, we have
\begin{equation}
\begin{split}
\int_{Q_+} P_1 z \partial_t z \,dxdt
&=
\int_{Q_+} \big(\partial_t^2 z - \Delta z + s^2\lambda^2 \varphi^2 z \big(|\partial_t \psi|^2 - |\nabla \psi|^2\big)\big) \partial_t z \,dxdt\\
&=\frac{1}{2}\int_{t=T}|\partial_tz|^2dS-\frac{1}{2\sqrt{2}}\int_{t=x\cdot \omega}|\partial_tz|^2dS-\int_{\Sigma_+}\nabla z\cdot\omega\partial_tzdS\\
&-\frac{1}{\sqrt{2}}\int_{t=x\cdot \omega}\nabla z\cdot\omega\partial_tzdS+\frac{1}{2}\int_{t=T}|\nabla z|^2dS-\frac{1}{2\sqrt{2}}\int_{t=x\cdot \omega}|\nabla z|^2dS\\
&+\frac{s^2\lambda^2}{2}\int_{t=T} \varphi^2\big(|\partial_t \psi|^2 - |\nabla \psi|^2\big) |z|^2 \,dS-\frac{s^2\lambda^2}{2\sqrt{2}}\int_{t=x\cdot \omega} \varphi^2\big(|\partial_t \psi|^2 - |\nabla \psi|^2\big) |z|^2 \,dS,
\end{split}
\end{equation}
which implies that
\begin{equation}
\begin{split}
&\frac{1}{\sqrt{2}}\int_{t=x\cdot\omega} |\partial_tz+\nabla z\cdot\omega|^2  dS +\frac{s^2\lambda^2}{\sqrt{2}}\int_{t=x\cdot \omega} \varphi^2\big(|\partial_t \psi|^2 - |\nabla \psi|^2\big) |z|^2 \,dS\\
&\leq \frac{1}{\sqrt{2}}\int_{t=x\cdot \omega}|\partial_tz|^2dS+\frac{2}{\sqrt{2}}\int_{t=x\cdot \omega}\nabla z\cdot\omega\partial_tzdS+\frac{1}{\sqrt{2}}\int_{t=x\cdot \omega}|\nabla z|^2dS\\
&
+\frac{s^2\lambda^2}{\sqrt{2}}\int_{t=x\cdot \omega} \varphi^2\big(|\partial_t \psi|^2 - |\nabla \psi|^2\big) |z|^2 \,dS\\
&=-2\int_{Q_+} P_1 z \partial_t z \,dxdt+\int_{t=T}|\partial_tz|^2dS+\int_{t=T}|\nabla z|^2dS\\
&+s^2\lambda^2\int_{t=T} \varphi^2\big(|\partial_t \psi|^2 - |\nabla \psi|^2\big) |z|^2 \,dS-s^2\lambda^2\int_{t=T} \partial_t\Big(\varphi^2\big(|\partial_t \psi|^2 - |\nabla \psi|^2\big)\Big) |z|^2 \,dS\\
&
-2\int_{\Sigma_+}\nabla z\cdot\omega\partial_tzdS.
\end{split}
\end{equation}
Using Cauchy-schwarz inequality, we obtain
\begin{equation}
\begin{split}
&\frac{s^{1/2}}{\sqrt{2}}\int_{t=x\cdot\omega} |\partial_tz+\nabla z\cdot\omega|^2  dS +\frac{s^{2/5}\lambda^2}{\sqrt{2}}\int_{t=x\cdot \omega} \varphi^2\big(|\partial_t \psi|^2 - |\nabla \psi|^2\big) |z|^2 \,dS\\
&\leq \int_{Q_+} |P_1 z|^2 \,dxdt+s^2\int_{Q_+} |\partial_t z|^2 \,dxdt+s^{1/2}\int_{t=T}|\partial_tz|^2dS+s^{1/2}\int_{t=T}|\nabla z|^2dS\\
&+s^{2/5}\lambda^2\int_{t=T} \varphi^2\big(|\partial_t \psi|^2 - |\nabla \psi|^2\big) |z|^2 \,dS-2\int_{\Sigma_+}\nabla z\cdot\omega\partial_tzdS.
\end{split}
\end{equation}
Combining (\ref{P12}) and \(\partial_tz=0\) on \(\Sigma_+\), we obtain
\begin{equation}\label{Pz}
\begin{split}
&s^{1/2}\int_{t=x\cdot\omega} |\partial_tz+\nabla z\cdot\omega|^2  dS +s^{2/5}\lambda^2\int_{t=x\cdot \omega} \varphi^2\big(|\partial_t \psi|^2 - |\nabla \psi|^2\big) |z|^2 \,dS\\
&
\leq  s\lambda \int_{Q_+} \varphi |\partial_t z|^2 \,dxdt
+  s\lambda \int_{Q_+} \varphi |\nabla z|^2 \,dxdt+ s^3\lambda^3 \int_{Q_+} F_{\lambda,t_0}(\phi)\varphi^3 | z|^2 \,dxdt\\
&\leq M\Big(\int_{Q_+}|Pz|^2dxdt+X_T+ s^{1/2}\int_{t=T}|\partial_tz|^2dS+s^{1/2}\int_{t=T}|\nabla z|^2dS+s^{5/2}\lambda^2\int_{t=T} \varphi^2\big(|\partial_t \psi|^2 - |\nabla \psi|^2\big) |z|^2 \,dS \Big).
\end{split}
\end{equation}
The expression for $X_{T}$ is
\begin{equation}\label{tT}
\begin{split}
X_{T}&=(\alpha - 1)s\lambda\int_{t=T}\partial_t z z\varphi (\partial_t^2\psi - \Delta\psi)dS-\frac{(\alpha - 1)s\lambda}{2}\int_{t=T} |z|^2\partial_t\varphi (\partial_t^2\psi - \Delta\psi)dS
\\
&-s\lambda^2\int_{t=T} \partial_t z  z\varphi(|\partial_t\psi|^2 - |\nabla\psi|^2)dS+\frac{s\lambda^2}{2}\int_{t=T} |z|^2 \partial_t\left(\varphi(|\partial_t\psi|^2 - |\nabla\psi|^2)\right)dS\\
&
-s\lambda \int_{t=T} |\partial_tz|^2\varphi\partial_t\psi dS
+2s\lambda \int_{t=T} \partial_tz \varphi\nabla z\cdot\nabla\psi dS
\\
&
-s\lambda \int_{t=T}|\nabla z|^2\varphi\partial_t\psi dS
-s^3\lambda^3\int_{t=T} |z|^2\varphi^3 (|\partial_t\psi|^2 - |\nabla\psi|^2)\partial_t\psi dS.
\end{split}
\end{equation}
When the time \(T\) is large enough, we claim that the terms at \(\{t=T\}\) can be removed of (\ref{Pz}). The proof is similar to the theorem in \cite{BBE}; we omit the details.

Since $z = v e^{s\varphi}$, we have
\[
\partial_t v+\nabla v\cdot\omega=e^{-s\varphi}(-s\lambda)\varphi (\partial_t \psi+\nabla v\cdot\psi)z+e^{-s\varphi}(\partial_t z+\nabla z\cdot\omega), \quad \text{in } Q_+.
\]
Hence 
\begin{equation}
\begin{split}
\int_{t=x\cdot\omega} e^{2s\varphi}|\partial_tv+\nabla v\cdot\omega|^2  dS \leq \int_{t=x\cdot\omega} |\partial_tz+\nabla z\cdot\omega|^2  dS+s^2\lambda^2\int_{t=x\cdot\omega} \varphi^2|\partial_t\psi+\nabla \psi\cdot\omega|^2 |z|^2 dS.
\end{split}
\end{equation}
Similarly, we have
\begin{equation}
e^{2s\varphi} |\partial_t v|^2 \le 2 |\partial_t z|^2 + 2 s^2 |\partial_t \varphi|^2 |z|^2
\le 2 |\partial_t z|^2 + 2 s^2 \lambda^2 \varphi^2|\partial_t\psi|^2 |z|^2 \quad
 \text{in } Q_+
\end{equation}
and
\begin{equation}
e^{2s\varphi} |\nabla v|^2 \le 2 |\nabla z|^2 + 2 s^2 |\nabla \varphi|^2 |z|^2
\le 2 |\nabla z|^2 + 2 s^2 \lambda^2 \varphi^2 |\nabla\psi|^2|z|^2\quad
 \text{in } Q_+.
\end{equation}
Therefore, inserting \(z = v e^{s\varphi}\) into (\ref{Pz}), we get
\begin{equation}\label{car}
\begin{split}
&s^{1/2}\int_{t=x\cdot\omega} e^{2s\varphi}|\partial_tv+\nabla v\cdot\omega|^2  dS +s^{2/5}\lambda^2\int_{t=x\cdot\omega} e^{2s\varphi }(|\partial_t\psi|^2-|\nabla\psi|^2)|v|^2dS\\
& + s  \int_{Q_+} e^{2s\varphi} \left( |\partial_t v|^2 + |\nabla v|^2 \right) dxdt + s^3  \int_{Q_+} e^{2s\varphi} |v|^2 dxdt \\
&
\leq M \int_{Q_+}e^{2s\varphi}|\partial_t^2v-\Delta v+Vv|^2dxdt.
\end{split}
\end{equation}
Now we take \(T\to\infty\) in (\ref{car}), which yields (\ref{Car1}).

\section{Proof of Theorem \ref{T2}}\label{Proof}

Let \(v=u_1-u_2\) and \(V=V_2-V_1\). Then, by proposition \ref{pro}, we have
\begin{equation}\label{v}
\partial_t^2 v -\Delta v+ V_1v = Vu_2, \quad (x, t) \in \mathbb{R}^n \times \mathbb{R}, \, t > x \cdot \omega, 
\end{equation}
\begin{equation}\label{v1}
v(x, x \cdot \omega, \omega) = -\frac{1}{2} \int_{-\infty}^0 V(x + \sigma \omega) d\sigma, \quad x \in \mathbb{R}^n, 
\end{equation}
\begin{equation}\label{v2}
V(x, t, \omega) = 0, \quad x \in \mathbb{R}^n, \, x \cdot \omega < t \ll 0. 
\end{equation}
On the plane \(t=x\cdot\omega\), The  relation (\ref{v1}) 
is equivalent to
\begin{equation}
\partial_tv+\nabla v\cdot\omega =-\frac{1}{2}V \quad \text{on the plane} \quad t=x\cdot\omega.
\end{equation}
which can be found in \cite{RU14}.
Applying (\ref{Car1}), we have 
\begin{equation}\label{car2}
s^{1/2}\int_{t=x\cdot\omega} e^{2s\varphi}|V|^2  dS \leq M \int_{Q_+^\infty}e^{2s\varphi}|V|^2dxdt.
\end{equation}
Taking the limit as \(t_0\to\infty\) in (\ref{car2}), we get
\begin{equation}
s^{1/2}\int_{t=x\cdot\omega}|V|^2  dS \leq M \int_{Q_+^\infty}|V|^2dxdt.
\end{equation}
This is
\begin{equation}
s^{1/2}\int_{B} |V|^2  dS \leq M \int_{B}|V|^2dx.
\end{equation}
Choosing \(s\) sufficiently large, we have
\begin{equation}
\int_{B}|V|^2  dx =0
\end{equation}
Then 
\[V=0 \quad \text{for all }x\in B.\]

\section*{Acknowledgment}
The work is supported by  the Shandong Provincial Natural Science Foundation (No. ZR2022QA111).

\end{document}